\overfullrule=0in
\def\St{{\hbox{\bf S}}^2}
\parskip = .05in \parindent = .2in
\def\v{\par\vskip .2in}  \def\n{\noindent}
\phantom \nobreak \v
\centerline{\bf The homotopy dimension of codiscrete subsets of the 2-sphere $\St$}
\v

\centerline{J. W. Cannon
\footnote*{This research was supported by NSF research grant DMS-10104030.}
and G. R. Conner}

\catcode`\"=0


\def\Rt{{\hbox{\bf R}}^2}

\def\Rn{{\hbox{\bf R}}^n}

\def\tH{\tilde{H}}
\def\inv{^{-1}}

\def"A{"mycala}
\def"B{"mycalb}
\def"C{"mycalc}
\def"D{"mycald}
\def"G{"mycalg}
"def"H{"hbox{"bf Hull}}

\def\inv{^{-1}}

\def\S{{\hbox{\bf S}}}

\def\Sn{{\hbox{\bf S}}^n}

\def\Bt{{\hbox{\bf B}}^2}

\def\Rt{{\hbox{\bf R}}^2}

\def\Rn{{\hbox{\bf R}}^n}

\def\mapright#1{\smash{\mathop{\longrightarrow}\limits^{#1}}}
\def\face(#1){\mapright{#1}}
"def"nt{"hbox{int}}
"def"cl{"hbox{cl}}
"def"e{"epsilon}
"def"diam{"hbox{diam}}
"v
{"narrower"narrower
"n{\bf Abstract.} Andreas Zastrow conjectured, and
Cannon-Conner-Zastrow proved, (see [3,pp. 44-45]) 
that filling one hole in the
Sierpinski curve with a disk results in a planar Peano
continuum that is not homotopy equivalent to a 1-dimensional
set. Zastrow's example is the motivation for this paper, where
we characterize those planar Peano continua that are
homotopy equivalent to 1-dimensional sets.

While many planar Peano continua are not homotopically 1-dimensional,
we prove that each has fundamental group that embeds in the
fundamental group of a 1-dimensional planar Peano continuum.

We leave open the following question: Is a planar Peano continuum
homotopically 1-dimensional if its fundamental group is isomorphic
with the fundamental group of a 1-dimensional planar Peano
continuum?

}
"v
"n{"bf 1. Introduction.}
We say that a subset $X$ of the 2-sphere $"St$ is {"sl codiscrete"/}
if and only if its complement $D(X)$, as subspace of $"St$, is discrete.
The set $B(X)$ of limit points of $D(X)$ in $"St$, which is necessarily
a closed subset of $X$ having dimension $"le 1$, is called the
{"sl bad set"/} of $X$. Our main theorem characterizes 
the homotopy dimension of
$X$ in terms of the interplay between $D(X)$ and $B(X)$:

{"bf Characterization Theorem 1.1.} 
Suppose that $X$ is a codiscrete subset of
the 2-sphere $"St$. Then $X$ is homotopically at most 1-dimensional if
and only if the following two conditions are satisfied.

(1) Every component of $"St "setminus B(X)$ contains a point
of $D(X)$.

(2) If $D$ is any closed disk in the 2-sphere $"St$, then the
components of $D "setminus B(X)$ that do not contain any point
of $D(X)$ form a null sequence.

[Recall that a sequence $C_1,C_2,"ldots$ is a {"sl null sequence"/}
if the diameters of the sets $C_n$ approach 0 as $n$ approaches
$"infty$.] Examples appear in Figures~1 and 2. Figure~1 gives
two examples of possible bad sets that are locally connected. The
one is a circle with countably many copies of the Hawaiian
earring attached. The other is a Sierpinski curve. The associated
codiscrete set will be homotopically 1-dimensional if and only
if condition~(1) is satisfied. Figure~2 gives an example of
a possible bad set that is not locally connected. In order that
the associated codiscrete set be homotopically 1-dimensional,
both conditions~(1) and (2) must be satisfied. Thus there must
be points of the discrete set near each point of the bad set
on both local sides of the bad set near the vertical limiting
arc.

\input epsf

\midinsert\epsfxsize=4in
\centerline{\epsfbox{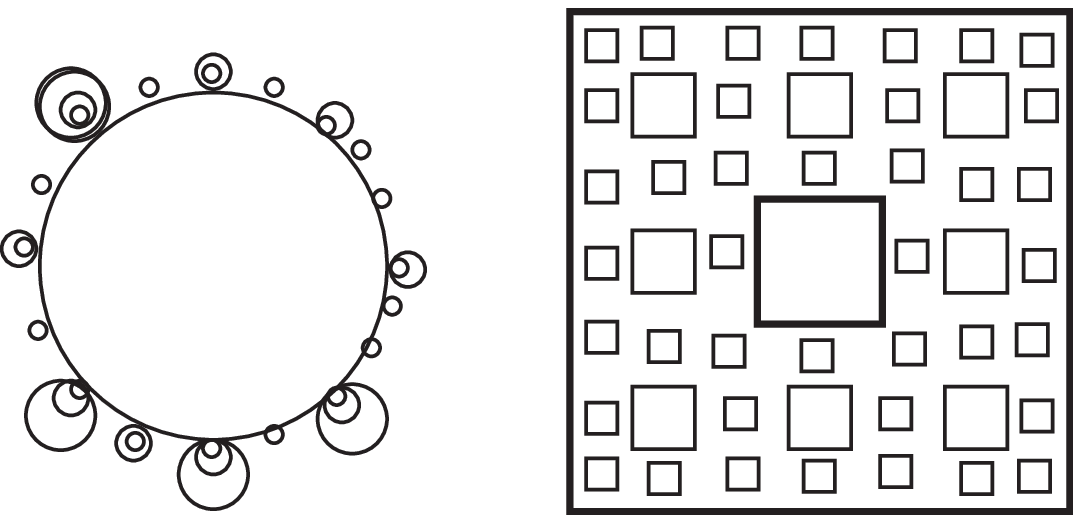}}
\centerline{Figure 1. Possible bad sets that are locally connected.}
\endinsert

\midinsert\epsfxsize=4in
\centerline{\epsfbox{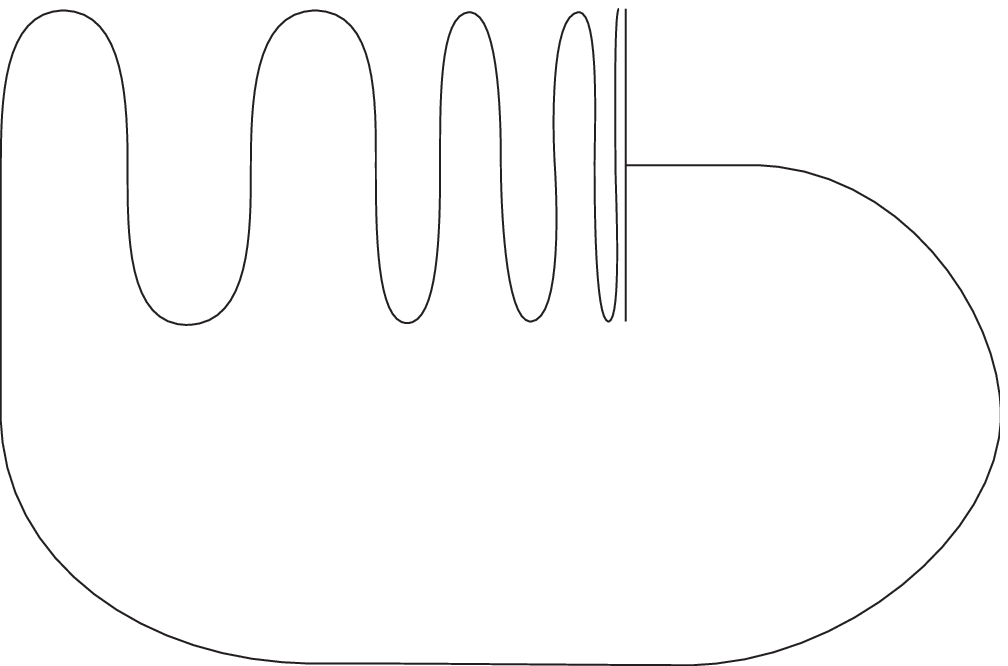}}
\centerline{Figure 2. A possible bad set that is not locally connected.}
\endinsert

A {"sl continuum} is a compact, connected metric space.
A {"sl Peano continuum} is a locally connected continuum;
equivalently, a Peano continuum is the metric continuous
image of the interval $[0,1]$.
Characterization Theorem~1.1  applies to all Peano continua in
the 2-sphere $"St$ because of the following well-known theorem:

{"bf Theorem 1.2.} Every Peano continuum $M$ in the 2-sphere $"St$
is homotopy equivalent to a codiscrete subset $X$ of $"St$.
Conversely, every codiscrete subset $X$ of $"St$ is homotopy
equivalent to a Peano continuum $M$ in $"St$.

We shall indicate later [after Theorem~2.4.2]
how this well-known theorem is proved.
For the moment, we simply mention that,
given $M$, one can obtain an appropriate codiscrete subset $X$
by choosing for $D(X)$ exactly one point from each component
of $"St "setminus M$. One can define the bad set $B(M)$ of $M$ as
the bad set $B(X)$ of $X$. It is natural to ask how restricted
bad sets are. The following theorem, which characterizes the
possible bad sets of codiscrete sets $X$, is actually an easy
exercise which we leave to the reader:

{"bf Theorem 1.3.} A subset $B$ of the 2-sphere $"St$ is the bad
set $B(X)$ of some codiscrete subset $X "subset "St$ if and
only if $B$ is closed and has dimension less than 2.

It is an easy matter to use Characterization Theorem 1.1 and
the construction inherent in Theorem 1.3 to construct all manner
of interesting planar Peano continua that are, or are not, homotopy
equivalent to a 1-dimensional set. All the examples that have
appeared in the literature (see [3] and [8]) are likewise easily checked by
means of Characterization Theorem 1.1.

In light of the fact that so many planar Peano continua are not
homotopically 1-dimensional, it is a little surprising to find
that their fundamental groups are essentially 1-dimensional
in the following sense:

{"bf Theorem 1.4.} If $M$ is a planar Peano continuum, then the
fundamental group of $M$ embeds in the fundamental group of a
1-dimensional planar Peano continuum. 

{"bf Corollary 1.4.1.} If $M$ is a planar Peano continuum, then
the fundamental group of $M$ embeds in an inverse limit of
finitely generated free groups.

{"bf Question 1.4.2.} If $M$ is a planar Peano continuum whose
fundamental group is isomorphic with the fundamental group of
some 1-dimensional planar Peano continuum, is it true that
$M$ is homotopically 1-dimensional?

The remaining sections of this paper
will be devoted to proofs of these theorems.

"n{\bf 2. Fundamental ideas and tools.} We collect here the basic
ideas and tools that will be used often in the proofs. Many
of these will be familiar to some of our readers. The topics
will be outlined in bold type so that the reader can quickly
find those topics with which they are not familiar. For many,
the best way to read the paper will be to turn immediately to
the later sections and return to this section only when they
encounter a tool or idea with which they are not familiar.

Our first fundamental idea is that {"bf homotopies of $X$ within itself
must fix the bad set $B(X)$ pointwise}. This general principle can
be applied to all connected planar sets $X$ and not just to codiscrete
sets. If $X$ is any connected planar set, then we may
define the {"sl bad set"/} $B(X)$ of
$X$ to be the set of points $x"in X$ having the property that, in
each neighborhood of $x$ there is a simple closed curve $J$ in $X$
such that the interior of $J$ in the plane $"Rt$ is not entirely
contained in the set $X$.

{"bf Theorem 2.1.} Suppose that $X$ is a connected planar set
and that $x"in B(X)$. Then every homotopy of $X$ within $X$
fixes the point $x$.

{"bf Proof.} Suppose that there is a homotopy $H:X "times [0,1] "to X$
such that $"forall y "in X$, $H(y,0) = y$ and such that $H(x,1) "ne x$.
Let $N_0$ and $N_1$ be disjoint neighborhoods of $x$ and $H(x,1)$,
respectively. By continuity, there is a neighborhood $M$ of $x$ in
$N_0$ such that $H(M,1) "subset N_1$. There is a round circle $J$ around
$x$ that is not contained in $X$ but intersects $X$ only in $N_0$.
There is a simple closed curve $K$ in $"nt(J) "cap M"subset X$ 
whose interior is not contained entirely in $X$. The annulus
$H(K "times [0,1])$ has its boundary components separated
by some component of $H"inv(J "cap X)$. That component maps into
a single component of $J"cap X$, where it can be filled in via
the Tietze Extension Theorem. This allows one to shrink $K$ in $X$,
an impossibility.

Our second fundamental idea is that of the {"bf convergence of a sequence
of sets}. Suppose that $A_1,A_2,"ldots$ is a sequence of subsets
of a space $S$. We say that a point $x"in S$ is an
element of $\liminf_i(A_i)$ if every neighborhood of $x$ intersects
all but finitely many of the sets $A_i$. We say that $x$ is an
element of $\limsup_i(A_i)$ if every neighborhood of $x$ intersects
infinitely many of the sets $A_i$. We say that the sequence $A_i$
converges if the $\liminf$ and $\limsup$ coincide. The limit is
defined to be this common $\liminf$ and $\limsup$. Here are the
fundamental facts about convergence, all of them well-known:

{"bf Theorem 2.2.1.} If $A_1, A_2,"ldots$ is any sequence of sets
in a separable metric space $S$, then there is a convergent
subsequence.

{"bf Proof.} Let $U_1,U_2,"ldots$ be a countable
basis for the topology of $S$. Let $S_0$ be the given sequence
$A_1,A_2,"ldots$ of subsets of the space $S$. Assume inductively
that a subsequence $S_i$ of $S$ has been chosen. If there
is a subsequence of $S_i$ no element of which intersects $U_{i+1}$,
let $S_{i+1}$ be such a subsequence. Otherwise, let $S_{i+1} = S_i$.
Let $S_"infty$ be the diagonal sequence, which takes as first element
the first element of $S_1$, as second element the second element
of $S_2$, etc. We claim that the subsequence $S_"infty$ of $S_0$
converges. Indeed, suppose that $x "in "limsup{S_"infty}$
That is, every neighborhood of $x$ intersects infinitely
many elements of $S_"infty$. Suppose that there is a neighborhood
$U_j$ of $x$ that misses infinitely many elements of $S_"infty$.
Then $S_j$, by definition, must miss $U_j$. But this implies
that all elements of $S_"infty$ with index as high as $j$ miss
$U_j$, a contradiction. Thus, every element of the $"limsup$
lies in the $"liminf$. Since the opposite inclusion is obvious,
these two limits are equal; and the sequence $S_"infty$ converges.

{"bf Theorem 2.2.2. Properties of the limit of a convergent sequence.} 
Suppose that the sequence $A_1,A_2,"ldots$ of nonempty subsets of a
separable metric space $S$ converges to a set $A$. Then,

(1) the set $A$ is closed in $S$;

(2) if $S$ is compact, then $A$ is nonempty and compact;

(3) If $S$ is compact and if each $A_i$ is connected, then
the limit $A$ is nonempty, compact, and connected.

(4) If $S$ is compact and if each $A_i$ has diameter $"ge "epsilon$,
then $A$ has diameter $"ge "epsilon$.

{"bf Proof.} Easy exercise.

We shall in more than one place make use of {"bf R. L. Moore's
Decomposition Theorem}.
In 1919 [9], R. L. Moore characterized
the Euclidean plane topologically. In 1925 [10], he noted that his axioms
were also satisfied by a large class of quotient spaces of the
plane, so that those identification spaces were also planes.

Since Moore's theorem is somewhat inaccessible to today's readers
because of evolving terminology and background, we will give a
fairly straightforward statement and outline the proof of this
theorem.

{\bf Moore Decomposition Theorem 2.3.1.} Suppose that $f:\St \to X$ is
a continuous map from the $2$-sphere $\St$ onto a Hausdorff space
$X$ such that, for each $x\in X$, the set $\St \setminus f\inv(x)$
is homeomorphic with the plane $\Rt$. Then $X$ is a $2$-sphere.

{\bf Remarks.} (1) The requirement that $\St \setminus f\inv(x)$ be
homeomorphic with $\Rt$ is equivalent to the requirement that
both $f\inv(x)$ and $\St \setminus f\inv(x)$ be nonempty and
connected.

(2) The theorem has the following generalization to higher dimensions:
{\sl Suppose that $f:\Sn \to X$ is a continuous map from the
$n$-sphere $\Sn$ onto a Hausdorff space $X$ such that, for each
$x\in X$, the set $\Sn \setminus f\inv(x)$ is homeomorphic with
the Euclidean space $\Rn$. Then $X$ is an $n$-sphere provided that,
in addition, $n\ge 5$, and $X$ satisfies the condition that
maps $g:\Bt \to X$ from the $2$-dimensional disk $\Bt$ into
$X$ can be approximated by embeddings.\/} This generalization
was conjectured and proved in many special cases by Cannon (see [1]
for a substantial discussion of these matters)
and proved in general by R. D. Edwards (see Daverman's book
[5].) The situation in dimensions $3$ and $4$ has
not been completely resolved.

The proof we shall give relies on a more intuitive theorem,
called the Zippin Characterization Theorem. (See, for example,
[13, p. 88].)

{\bf Zippin Characterization Theorem 2.3.2.}The space $X$ is a $2$-sphere
if the following four conditions are satisfied:

(i) $X$ is a nondegenerate Peano continuum.

(ii) No point $x\in X$ separates $X$ (so that, in particular,
$X$ contains at least one simple closed curve).

(iii) Each simple closed curve $J\subset X$ separates $X$.

(iv) No arc $A\subset X$ separates $X$.

{\bf Proof of the Moore Decomposition Theorem on the basis
of the Zippin Characterization Theorem.} We verify the four
conditions of the Zippin Theorem in turn. (Note that conditions
(iii) and (iv) are true in the $2$-sphere by standard
homological arguments. We shall use those same arguments here.)

(i): Since $X$ is Hausdorff, the map $f$ is a closed surjection;
hence it is easy to verify the conditions of the Urysohn metrization
theorem so that $X$ is metric. (See [11,
Theorem 34.1].) Since $\St$ is a Peano continuum, that is,
a metric continuous image of $[0,1]$,
so also is $X$. Since, $\forall x \in X$,
both $f\inv(x)$ and $\St \setminus f\inv(x)$ are nonempty, $X$ has more
than one point; that is, $X$ is nondegenerate.

(ii): By hypothesis, $\St \setminus f\inv(x)$ is connected. Hence
$X \setminus \{x\} = f(\St\setminus f\inv(x))$ is also connected.

(iii): Let $p_1,p_2\in J$ cut $J$ into two arcs $A_1$ and $A_2$. Then
$f\inv(A_1)$ and $f\inv(A_2)$ are compact, connected, and have
nonconnected intersection $f\inv(p_1)\cup f\inv(p_2)$. The
reduced Mayer-Vietoris homology sequence for the pair
$U = \St \setminus f\inv(A_1)$ and $V = \St \setminus f\inv(A_2)$
contains the segment
$$H_1(\St \setminus f\inv(A_1)) \oplus H_1(\St \setminus f\inv(A_2)) 
\to H_1(\St \setminus (f\inv(p_1)\cup f\inv(p_2)) \to \tH_0(\St \setminus (f\inv(J)),$$
where $H_1(U) = H_1(V) = 0$ since $f\inv(A_1)$ and $f\inv(A_2)$ are
connected and $H_1(U\cup V) \ne 0$ since $f\inv(A_1) \cap f\inv(A_2)$
is not connected. Thus $\tH_0(\St\setminus f\inv(J)) =\tH_0(U\cap V) \ne 0$,
so that $f\inv(J)$ separates $\St$. Consequently, $J$ separates $X$.

(iv): If $p\in A$ separates $A$ into arcs $A_1$ and $A_2$, and if $A$
separates $x$ and $y$ in $X$, then we claim that one of $A_1$ and
$A_2$ also separates $x$ and $y$ in $X$; indeed, we see this by considering
$f\inv(A) = f\inv(A_1) \cup f\inv(A_2)$, which must separate $f\inv(x)$
from $f\inv(y)$ in $\St$. The
reduced Mayer-Vietoris homology sequence for the pair
$U = \St \setminus f\inv(A_1)$ and $V = \St \setminus f\inv(A_2)$
contains the segment
$$0 \to \tH_0(\St\setminus f\inv(A))\to \tH_0(\St\setminus A_1)\oplus \tH_0(\St\setminus A_2).$$
The element $x-y$ represents a nonzero element of the center group,
hence maps to a nonzero element of $\tH_0(\St\setminus A_1)\oplus\tH_0(\St\setminus A_2)$,
as desired.

By induction, one obtains intervals $I_0\supset I_1\supset\cdots$
that separate $x$ and $y$ in $X$ such that $\cap_{n=1}^\infty I_n$
is a single point $q$ that does not separate $x$ from $y$. But an
arc $\alpha$ from $x$ to $y$ in the path connected open set $X\setminus \{q\}$
misses some $I_n$, a contradiction. We conclude that $A$ cannot
separate $X$.

The proof of the Moore Decomposition Theorem 2.3.1 is complete.

Our fourth topic is that of {"bf locally connected continua in the
plane}. 

{"bf Theorem 2.4.1.} Suppose that $M$ is a continuum (= compact,
connected subset) in the 2-sphere $"St$. Then $M$ is a Peano
continuum (= locally connected continuum) if and only if the
following four equivalent conditions are satisfied.

(1) For each disk $D$ in $"St$, the components of $D "setminus M$
form a null sequence.

(1$'$) For each disk $D$ in $"St$, the components of $D "cap M$
form a null sequence.

(2) For each annulus $A$ in $"St$, the components of $A "setminus M$
that intersect both boundary components of $A$ are finite
in number.

(2$'$) For each annulus $A$ in $"St$, the components of $A "cap M$
that intersect both boundary components of $A$ are finite in
number.

{"bf Proof.} Assume that $M$ is locally connected but that
(1) is not satisfied, so that, for some
disk $D$ in $"St$, the components of $D"setminus M$ do not form
a null sequence. Then some sequence $U_i$ of such components
converges to a nondegenerate continuum $U$ in $"St$ by
Theorems 2.2.1 and 2.2.2. Let $A$ be an annulus in $"St$
that separates two points of $U$. Then each $U_i$ contains
an arc $A_i$ irreducibly joining the two ends of $A$.
We may assume that they converge to a continuum $A'$ joining the
two ends of A. The continuum $A'$ must be a subset of $M$, for
otherwise it could not have points of infinitely many of the
components $U_i$ close to it. Since the arcs $A_i$ converge
to $A'$, there must be two of them, which we may
number as $A_1$ and $A_2$, that have no other $A_i$ nor
$A'$ between them. There are then only two components of $A "setminus
(A' "cup A_1 "cup A_2)$ that can contain any of the remaining
$A_i$. This allows us to choose a subsequence, which we may assume
is the sequence $A_3,A_4,"ldots$, such that each $A_i$ is
adjacent to $A_{i+1}$, with neither $A'$ nor any other $A_j$
between them. They must therefore be separated by a
component $M_i$ of $A "cap M$ that intersects both ends of $A$.
The components $M_i$ converge to a subcontinuum of $A'$ that
joins the ends of $A$. This shows that $M$ is not locally
connected at these points of $A'$, a contradiction.

Suppose (1) is satisfied but (1$'$) is not. That is, there is a disk
$D$ in $"St$ and infinitely many large components of $D "cap M$.
We may take a sequence of such components that converge
to a nondegenerate subcontinuum of $M$.
We take an annulus $A$ that separates two points of the limit continuum.
Infinitely many of the large components cross this annulus.
They are separated by large components of $A"setminus M$ that
cross the annulus. Arcs in these components that cross the 
annulus allow one to form a disk $D$ that is crossed by infinitely
many large components of $D "setminus M$, a contradiction to (1).
We conclude that (1$'$) is satisfied.

Similar arguments show that (1$'$) implies (1) and that these
are equivalent to (2) and (2$'$).

Finally, if $M$ is not locally connected, then there is a point
$p"in M$ and a neighborhood $N$ of $p$ in $M$ such that
$p$ is a limit point of the components of $N "cap M$ that do
not contain $p$. Each of these components intersects the boundary
of $M$. These large components contradict (1$'$).

{"bf Theorem 2.4.2.} Suppose that $M$ is a Peano continuum in the
2-sphere $"St$, and suppose that $U$ is a component of the complement
of $M$ in $"St$. Then there is a map $f:"Bt "to "cl(U)$
from the 2-disk $"Bt$ onto the closure of the domain $U$ that
takes $"nt(B^2)$ homeomorphically onto $U$ and takes
$"So = "bd("Bt)$ continuously onto $"bd(U)$. In addition, if $A$
is a free boundary arc of $"cl(U)$, then we may assume that
the map $f$ is one to one over the arc $A$.

{"bf Remark.}
That the arc $A$ is {"sl free} means that $A$ is accessible
from precisely one of its sides from the domain $U$
and that $"nt(A)$ is an open subset of $"bd(U)$.

{"bf Indication of proof.} There are well-known, completely
topological proofs of this theorem. However, refinements of
the Riemann Mapping Theorem also give very enlightening
analytic information. The relevant analytic theory is the theory
of {"sl prime ends}. There is a good exposition of the theory
in John B. Conway's readily available textbook, [4,
Chapter~14, Sections 1-5].
It follows from the local connectivity
of $M$ (applying Theorem 2.4.1(1))
that the {"sl impressions} of the prime ends in $U$ are all 
singletons. By the
theory of prime ends, the Riemann mapping from $"nt("Bt)$
onto $U$ extends continuously to the boundary. 

That the arc $A$ is {"sl free} means that $A$ is accessible
from precisely one of its sides from the domain $U$
and that $"nt(A)$ is an open subset of $"bd(U)$.
Consequently, the prime ends
at $A$ correspond exactly to the points of $A$ so that
the map is one to one over $A$.

{"bf Proof of Theorem 1.2.} Suppose that $M$ is a locally
connected continuum in $"St$. If $M = "St$, then $M$ is
already codiscrete. Otherwise, let $U_1, U_2, "ldots$ denote
the complementary domains of $M$ in $"St$. By Theorem~2.4.1,
the components of $"St "setminus M$ form a null sequence.
By Theorem~2.4.2, there is for each $i$ a continuous
surjection $f_i:"Bt "to "cl(U_i)$ that takes $"So$ onto
the boundary of $U_i$ and takes the interior of $"Bt$ 
homeomorphically onto $U_i$. Let $p_i = f_i(0)$. Then
the set $D = "{p_1,p_2,"ldots"}$ is obviously discrete.
The set $"cl(U_i)"setminus "{p_i"}$ can obviously be
deformed into the boundary of $U_i$ by pushing points
away from $p_i$ along the images under $f_i$ of radii in $"Bt$.
These deformations can be combined to deform all
of $X = "St "setminus D$ onto $M$ since the $U_i$ form
a null sequence. We conclude that $M$ is homotopy equivalent
to the codiscrete set $X = "St "setminus D$.

Conversely, if $X$ is codiscrete, then we may take, about the
points $p$ of $D(X)$, small disjoint round disks $d(p)$.
The continuum $M = "St "setminus "cup_p"nt(d(p))$ is
a Peano continuum to which $X$ can be deformed by a 
strong deformation retraction.

This completes the proof of Theorem 1.2.

{"bf We may think of the proof of Characterization Theorem 1.1 as a
substantial generalization of the proof of Theorem 2.4.2. We shall need an
intermediate generalization of Theorem 2.4.2
that deals with compact sets that
act much like Peano continua but are not necessarily connected.
We shall deal with them by joining them together by arcs so
as to form a Peano continuum.}

{"bf Definition 2.5.1.} A connected open subset $U$ of $"St$
is called a {"sl Peano domain"/} if its nondegenerate boundary
components form a null sequence of Peano continua. [Note that
there may be uncountably many additional components that
are single points.]

{"bf Theorem 2.5.2.} Suppose that $U$ is a connected open subset
of the 2-sphere $"St$. Then the following three conditions
are equivalent:

(1) The open set $U$ is a Peano domain.

(2) For each disk $D$ in $"St$, the components of $U "cap D$ form
a null sequence.

(3) There is a continuous surjection $f:"Bt "to "cl(U)$
such that $f("So)"supset "bd(U)$ and $f|"nt("Bt)$ is
a homeomorphism onto its image.

{"bf Remark.} Note that (1) generalizes the
notion of local connectedness. Note that (2) generalizes
characterization (1) of local connectedness in Theorem
2.4.1; the reader
can reformulate (2) in each of the ways suggested
by Theorem 2.4.1. Note that (3) generalizes Theorem 2.4.2.  
Note that, in the proof, we can assume
that the map $f$ is 1-1 over given free boundary arcs
of $U$ because the same thing is true in Theorem 2.4.2.

{"bf Proof.} Assume (1), so that $U$ is a Peano domain. 
Assume that (2) is
not satisfied, so that there is a disk $D$ in $"St$
such that the components of $U "cap D$ do not form
a null sequence. Then some sequence $U_1,U_2,"ldots$ of
components converges to a nondegenerate continuum $M$. The
continuum $M$ must be a subset of a boundary component of
$U$. We may assume that the components $U_1,U_2,"ldots$
are separated from each other by large boundary components
of $U$. There are only finitely many large boundary components
of $U$. Hence infinitely many of the separators must come
from the same boundary component. It follows that the
limit, namely $M$, is also in the same boundary component.
But this boundary component is not locally connected at
the points of $M$, a contradiction. We conclude that (2) is
satisfied so that (1) implies (2).

Assume that (2) is satisfied. Assume that (1) is not satisfied.
Then either there is a component of $"bd(U)$ that is not
locally connected, or there exist infinitely many components
of $"bd(U)$ having diameter $"ge "epsilon$, for some 
fixed $"epsilon > 0$. In either case, taking a convergent
sequence of large components, we find the existence of an
annulus $A$ in $"St$ and components $X_1,X_2,"ldots$ of
$"bd(U) "cap A$, each of which intersects both components of $"bd(A)$.
These components $"bd(U) "cap A$ must be separated by
large components of $A"cap U$. if we remove a slice from one
of these large separating components, we obtain a disk $D$
that is crossed by infinitely many large components of $U "cap D$,
which contradicts (2). Therefore (2) implies (1).

Assume that (3) is satisfied, so that there is a continuous 
surjection $f:"Bt "to "cl(U)$
such that $f("So)"supset "bd(U)$ and $f|"nt("Bt)$ is
a homeomorphism onto its image. Assume that (1) is not satisfied,
so that there is either a component of $"bd(U)$ that is not
locally connected, or, there exist infinitely many
components of $"bd(U)$ each having diameter greater than
some fixed positive number $"epsilon$. In either case, we find
by taking limits that there is an annulus $A$ in $"St$ and
components $X_1,X_2,"ldots$ of $"bd(U) "cap A$, each of which intersects
both components of $"bd(A)$. We may assume that $X_1,X_2,"ldots$ 
converges to a continuum $X_0$ joining both components of $"bd(A)$.
We may assume that $X_{i-1} "cup X_{i+1}$ 
separates $X_i$ from $X_0$ in $A$, for $i = 2,3,"ldots$.

Pick $p_i"in X_i "cap "nt(A)$ such that $p_1,p_2,"ldots "to p_0$. 
Let $q_0,q_1,q_2,"ldots"in "So$ be points such that $f(q_i) = p_i$.
Let $B_i$ be the straight-line segment in $"Bt$ joining $q_0$ to
$q_i$. We may assume that the arcs $B_i$ converge to an arc
or point $B$ in $"Bt$. We shall obtain a contradiction as follows.

The image $f(B_i)$ joins $X_i$ to $X_0$. It
misses $X_{i-1}"cup X_{i+1}"subset "bd(U)$ 
since $f(q_i)"in X_i$, $f(q_0) "in X_0$, and $f("nt(B_i)
"subset U$. Hence, traversing $B_i$ from $q_i$ toward $q_0$,
there exists a first point $b_i "in B_i$ such that
$f(b_i)"in "bd(A)$. We may assume that $b_i"to b"in "Bt$
and $f(b_i)"to f(b_0)"in "bd(A)$. Since $f(b_i)$ is separated from
$X_0$ by $X_{i-1}"cup X_{i+1}$ in $A$ and since $X_i"to X_0$,
we may conclude that $f(b_0)"in X "cap "bd(A)$. Hence 
$b_0 "in "So"setminus"{q_0"}$. But $b_0$ must therefore be an endpoint
of $B$ distinct from $q_0$ and must therefore be the limit of
the points $q_i$. We find that $f(q_i)"to p_0"in "nt(A)$ and
$f(q_i)"to f(p_0) "in "bd(A)$, a contradiction.

We conclude that (3) implies (1).

It remains to prove that (1) implies (3). This is by far the
hardest of the implications. It is a generalization of 
the rather deep Theorem~2.4.2, and we shall reduce it to that
theorem. We shall also make use of the wonderful
R.L. Moore Decomposition Theorem~2.3.1.

Our plan is to connect $"bd(U)$ by deleting from $U$ a
null sequence $A_1,A_2, "ldots$ of arcs to form a new
connected open set $V= U "setminus "cup_iA_i$ whose
boundary $"bd(V)= "bd(U) "cup "bigcup_iA_i$ is a locally
connected continuum. Then we simply apply Theorem~2.4.2.

For convenience, we smooth the nondegenerate
components $C$ of $"bd(U)$ as follows. We define $U_C$ 
to be the component of $"St "setminus C$ that contains $U$.
Since $C$ is locally connected by (1), we may
apply Theorem~2.4.2 to find a continuous surjection 
$g:"Bt "to C "cup U_C$ that takes $"So$ onto $C$ and takes $"nt("Bt)$
homeomorphically onto $U_C$.  Thus, pulling 
$U_C$ radially into itself along the
images of radii, we find that we lose no generality in assuming
that $C$ is a topological circle. Since the nondegenerate components
of $"bd(U)$ form a null sequence by (1), we may repeat the
argument infinitely often to conclude that we lose no generality
in assuming that each nondegenerate component is a simple
closed curve. That is, $U$ is the complement of a null
sequence of disks $D_1,D_2,"ldots$ and a 0-dimensional set $D$,
the union of $D_1,D_2,"ldots$, and $D$ being closed. 

We wish to construct a nice sequence of cellulations
of the 2-sphere that respect the boundary components of
$U$. If, for example, we wish to concentrate on some
particular finite set $S$ of the large disks $D_i$,
we may form an upper semicontinuous decomposition
of $"St$ by declaring the other $D_i$'s that miss $S$ to be
the nondegenerate elements of the decomposition. By R.L. Moore's
Decomposition Theorem 2.3.1, the quotient space is the
2-sphere $"St$. The (homeomorphic) image of $U$ in this
new copy of $"St$ will have, as complement, the (images of the)
elements of $S$ and a 0-dimensional set that is closed
away from $S$. It is then an easy matter to cellulate $"St$
so that the elements of $S$ cover a subcomplex and
the remainder of the 1-skeleton misses $"bd U$ entirely.

As a consequence, we find that there is a sequence $S_1,S_2,"ldots$
of arbitrarily fine cellulations of $"St$, $S_{i+1}$ subdividing
$S_i$, such that, for each $i$, the following conditions
are satisfied:

(i) Two 2-cells of $S_i$ that intersect intersect in an arc.

(ii) The 1-skeleton of $S_i$ misses all of the 0-dimensional
part $D$ of $"bd(U)$.

(iii) $"forall j$, the 1-skeleton of $S_i$ either misses the
disk $D_j$ or contains $"bd(D_j)$. Consequently, $S_i$ has
a distinguished finite subcollection of disks $D_j$ that
are precisely equal to unions of 2-cells of $S_i$. All other
disks $D_k$ will lie in the interiors of 2-cells of $S_i$.

(iv) If a 2-cell $C$ of $S_i$ has a boundary point in some
$"bd(D_j)$, with $"nt(C) "not"subset D_j$, then
$"bd(C) "cap "bigcup_kD_k$ is an arc in $"bd(D_j)$.

We shall string the components of $"bd U$ together by arcs
that run through $U$. These arcs will be built by approximation.
The $i$th approximation will consist of arcs that join
certain 2-cells of the cellulation $S_i$. 

It is necessary to distinguish four types of 2-cells 
in the cellulation $S_i$:

A 2-cell $C$ of $S_i$ is of type~0 if it lies entirely
in $U$.

A 2-cell $C$ is of type~1 if 
it lies entirely in the complement of $U$, hence lies in
one of the distinguished disks $D_j$ of the cellulation
$S_i$ (see (iii) above). 

A 2-cell $C$ is of type~2 if it intersects both $U$ and
the complement of $U$, but its boundary lies entirely in $U$.

A 2-cell $C$ is of type~3 if its boundary intersects both $U$ and
the complement of $U$. Condition~(iv) above implies that a 2-cell
$C$ of type~3 has boundary that intersects precisely one
disk $D_j$, that $D_j$ is one of the distinguished disks of $S_i$,
and the intersection is a boundary arc of each.

We shall essentially ignore the 2-cells of type~0.
We shall deal with the disks of type~1 only implicitly
by considering instead their unions that give the distinguished disks
$D_j$ of the cellulation $S_i$ (see (iii) above). 
Cells of type~2 will be joined to these distinguished disks by
arcs in $U$. Cells of type~3 will be joined to these
distinguished disks by their intersecting boundary arcs.

It will be convenient to use the notation $C^*$ for the
union of the elements of a collection $C$ of sets.

Let $"D_1$ denote the collection of $D_j$'s that are
distinguished in the cellulation $S_1$. Then
$"D_1^* = "cup"{D"in "D_1"}$. We may assume
$D_1 "in "D_1$. We may pick a collection of arcs $"A_1$ from 
the 1-skeleton $S_1^{(1)}$ of $S_1$ that irreducibly
joins together these distinguished disks $D_j"in "D_1$
and the cells of $S_1$ of type~2. Then $"C_1 = "D_1^* "cup "A_1^*$
is a contractible set.

All of the cells of $S_1$ of type~1 are contained in $"C_1$.
We may consider all cells of type~2 and 3 as attached to this
contractible set in the following way. For each cell $C$ of
type~2, pick one point of intersection with $"C_1$ as attaching point.
For each cell $C$ of type~3, pick as attaching arc the
boundary arc of $C$ that lies in a distinguished disk.

We proceed by induction. We assume that we have constructed
contractible sets $"C_1"subset "C_2"subset "cdots "C_i$, that lie
except for distinguished disks of $S_1,S_2,"ldots,S_i$, in the
1-skeletons of the cellulations. We may impose one additional
condition on the cellulation $S_{i+1}$:

(v) For each cell $C$ of $S_i$ that has type~2 or 3, that part
of the 1-skeleton of $S_{i+1}$ that lies in the interior
of $C$, taken together with the attaching point (type~2) or
attaching arc (type~3), is connected.

All of the action in creating $"C_{i+1}$ takes place in
the individual cells $C$ of $S_i$ of type~2 and 3. 
We may pick a collection of arcs $"A_{i+1}(C)$ from 
that part of the 1-skeleton of $S_{i+1}$ that lies in the interior
of $C$, taken together with the attaching point (type~2) or
attaching arc (type~3), that irreducibly
joins together the attaching set of $C$,
the distinguished disks $D_j"in "D_{i+1}$ in $C$,
and the cells of $S_{i+1}$ of type~2 in $C$. All of these
new distinguished disks and all of these new arcs can be
added to $"C_{i}$ to form a new contractible set $"C_{i+1}$.
We denote the entire union $"bigcup_C"A_{i+1}(C)$
of arcs as $"A_{i+1}$.

For each of the new cells of types~2 and 3, we choose an
attaching point or arc as before.

We leave it to the reader to verify that 
$M = ("St "setminus U) "cup "bigcup_i(A_i)$ is a single
locally connected continuum with a single complementary
domain $V = U "setminus "bigcup_i(A_i)$. 

By Theorem 2.4.2, there is a map $f:"Bt "to "cl(V)$
from the 2-disk $"Bt$ onto the closure of the domain $V$ that
takes $"nt(B^2)$ homeomorphically onto $V$ and takes
$"So = "bd("Bt)$ continuously onto $"bd(V)$. The same
map establishes condition (3) of Theorem 2.5.2.

This completes the proof that (1) implies (3). Thus all three
conditions of Theorem~2.5.2 are equivalent, as claimed.
The proof of Theorem 2.5.2 is therefore complete.

Our final theorem of this section shows how to push a
Peano domain onto its boundary together with a 1-dimensional
set provided the domain is punctured on a nonempty discrete
set. This easy theorem will be needed as the last step in the
proof of Theorem~1.1.

{"bf Theorem 2.5.3.} Suppose that $U$ is a Peano domain in $"St$ and
that $C$ is a nonempty countable or finite subset of $U$ that
has no limit points in $U$. Then $"cl(U) "setminus C$ can
be retracted by a strong deformation retraction onto a
1-dimensional set that contains $"bd(U)$.

{"bf Proof.} By Theorem 2.5.2, we know that there is a
continuous surjection $f:"Bt "to "cl(U)$
such that $f("So)"supset "bd(U)$ and $f|"nt("Bt)$ is
a homeomorphism onto its image. 

Since $f("nt("Bt))$ is dense in $f("Bt) = "cl(U)$
and disjoint from $f("So)$, $f("So)$ must be 1-dimensional.
Hence it is an easy exercise to show that we may modify $f$
slightly over $U$ so that $f("So)$ misses $C$. We may further
modify $f$ so that $f$ maps the origin $0 "in "Bt$ to a point
of $C$ and so that all other points of $C$ have preimages on different
radii of $"Bt$. Let $f"inv(C) = "{c_0 = 0, c_1, c_2, c_3,"ldots"}$.
Let $A_1,A_2,"ldots$ be the radial arcs beginning at $c_1,c_2,"ldots$,
respectively, and ending on $"So = "bd("Bt)$. Let $D_1,D_2,"ldots$
be disjoint round disks in $"nt("Bt) "setminus "{0"}$
centered at $c_1,c_2,"ldots$, respectively, such that the only
$A_j$ intersected by $D_i$ is $A_i$. Let 
$V = "nt("Bt)"setminus [ "bigcup_i A_i "cup "bigcup_i D_i]$. 
Then $"Bt "setminus f"inv(C)$ can obviously be retracted by a strong
deformation retraction onto the 1-dimensional set $"bd(V)$.
Hence $f("Bt) "setminus C = "cl(U) "setminus C$ can be retracted
by a strong deformation retraction onto the 1-dimensional set $f("bd(V))$.

"n{\bf 3. The necessity of conditions (1) and (2) in Theorem 1.1.}

We assume that $X$ is a codiscrete set that is homotopy equivalent
to a metric 1-dimensional set $Y$. Let $f:X"to Y$ and $g:Y "to X$ be
homotopy inverses.

We isolate the three key technical constructions as lemmas.
Each of these is standard and well-known. We omit the proofs.

{"bf Dimension Lemma 3.1.} If $g:Z' "to Z$ is any map from a 
1-dimensional compactum $Z'$ into the closure $Z$ of an
open subset $U$ of $"St$, then $g$ is homotopic, by a homotopy
which only moves points in $U$ to a map $g':Z'"to Z$ such that
$g'(Z')"cap U$ is 1-dimensional. [The key ideas are explained,
for example, in [12, Exercises for Chapter 3, Sections G and H].]

{"bf Homotopy Lemma 3.2.} (i) Let $C "subset "St$ be closed,
and let $H:C "times [0,1] "to "St$ denote a deformation
of $C$ that begins at the identity (that is, $"forall c"in C$,
$H(c,0)=c$). Then $H$ can be extended to a deformation
$H':"St "times [0,1] "to "St$. (ii) If $H$ moves no
point as far as $"e > 0$, then we may require that $H'$ have the
same property. (iii) If $N$ is an open set containing the
support of $H|"bd C$ in $"St$, then we may require that $N$
contain the support of $H'|"St"setminus C$. [See [, \S 62, Lemma~62.1
and Exercise~3.]

{"bf Ring Lemma 3.3.} Suppose condition (2) of Theorem 1.1 fails.
Then there are a ring $R'$ in $"St$ and components 
$U_1',U_2',"ldots$ of $R' "setminus B(X)$ such that each
$U_j'$ intersects both boundary components of $R'$ and
misses the set $D(X)$. [See Theorem~2.4.1 and its proof.]

{"bf The three lemmas imply the theorem as follows:} By precomposing
the homotopy equivalence $f$ with a deformation retraction onto
a compact subset of $X$, we may assume that the
image $f(X)$ is a 1-dimensional continuum $Z'$.  
By Dimension Lemma 3.1, we may assume that
$g"circ f(X)"setminus B(X)$ is 1-dimensional. Let $G:X"times[0,1]"to "St$
be a homotopy that begins with the identity on $X$ and ends
with $g"circ f$. By Theorem~2.1, we see that $G(x,t)=x$
for each $x"in B(X)$.

Assume that condition (1) of the hypothesis of Theorem 1.1 fails, so that
there is some component $U$ of $"St "setminus B(X)$ contains no point
of $D(X)$. Hence $U "subset X$. Let $H:"cl(U) "times [0,1] "to "St$
denote the restriction of $G$ to $"cl(U)"times[0,1]$. Since
$H$ fixes $"bd U "subset B(X)$, we may extend
$H$ to a deformation $H'$ of $"St$ that fixes $"St "setminus U$
pointwise. Since 
$H'("St "times "{1"})"cap U "subset G("St "times"{1"})"cap U$
is 1-dimensional, we see that $H'$ deforms $"St$ into a proper
subset of itself, which is impossible. Hence condition (1) must
be satisfied.

Assume that condition (2) of the hypothesis of Theorem 1.1 fails.
Then, by Ring Lemma 3.3, there are a ring $R'$ in $"St$ and
components $U_1',U_2',"ldots$ of $R' "setminus B(X)$ such that
each $U_j'$ intersects both boundary components of $R'$ and
fails to intersect the set $D(X)$.

By passing to a subsequence, we may assume that the components 
$U_1',U_2',"ldots$ converge to a continuum $A$ that joins the
two boundary components of $R'$. Since the components $U_j'$ are
separated by $B(X)$, it follows that $A "subset B(X)$.
Let $D$ be a small disk in $"nt(R')$ centered at some
point of $A$. Since the deformation $G$ 
constructed above moves no point of $B(X)$, there is a 
neighborhood $N$ of $A$ in $X$, no point of which is moved by
$G$ as far as 1/2 the distance from $"bd R$ to $D$. We
shoose $j$ so large that $"cl(U_j)"subset N$ and
$U_j "cap "nt(D) "ne "emptyset$. Since no point of $D(X)$ lies
in $U_j$, all of $"cl(U_j)$ lies in $X$. 

We let $H:"cl(U_j)"times[0,1] "to "St$ be the restriction of
$G$ to $"cl(U_j)"times [0,1]$. By Homotopy Lemma~3.2(i), there
is a deformation $H':"St"times[0,1]"to "St$ that extends
$H$. By Homotopy Lemma~3.2(iii), we may require that
$H'|["St"setminus "cl(U_j)]"times [0,1]$ move points only
near $("bd R)"cap "cl(U_j)$, a set that contains the
support of $H|"bd U_j "times [0,1]$. By Homotopy Lemma~3.2(iii),
we may require that no points of $"St "setminus U_j$ be carried
into $D "cap U_j$. Hence $H'$ is a homotopy of $"St$ that takes
$"St$ to a proper subset of itself, an impossibility. Hence
condition (2) of Theorem~1.1 is also satisfied.

"n{\bf 4. The sufficiency of conditions (1) and (2) in 
Characterization Theorem~1.1.}

We assume conditions (1) and (2) of Characterization Theorem~1.1.
That is, the open set $U_0 = "St "setminus B(X)$ satisfies the
following two conditions:

(1) Each component of $U_0$ contains a point of $D(X)$.

(2) If $D$ is any disk in $"St$, then the components
of $U_0 "cap D$ that contain no point of $D(X)$ form
a null sequence.

Our goal is to show that $X$ is homotopy equivalent to a 1-dimensional
set.

Notice that properties (1) and (2) make no explicit mention of
the bad set $B(X)$ and are simply properties that an open subset
of $"St$ may or may not have. This is an important observation,
because our proof that $X$ is homotopy equivalent to a 1-dimensional
set will involve a complicated induction that will involve
a decreasing sequence $U_0 "supset U_1 "supset U_2 "supset "cdots$ of
open sets, each of which satisfies properties (1) and (2).

It will also be convenient to adopt the following terminology:
we say that set is {"sl punctured} if it contains
a point of $D(X)$. Otherwise, we say that it is {"sl unpunctured}.

We first have to deal with the trivial case where $B(X) = "emptyset$.
If $B(X) = "emptyset$, then the single 
component $"St = "St "setminus B(X)$
must contain a point of $D(X)$ by (1). Thus there must be
at least one point of $D(X)$ and at most finitely many. Hence $X$
is clearly homotopy equivalent to a point or bouquet of circles.

From now on, we may assume that the set $D(X)$ is infinite and
the set $B(X)$ is nonempty.
Since $D(X)$ is countable, we may list the points $p_0,p_1,p_2,"ldots$
of $D(X)$. We need to show that $X$ is 
homotopy equivalent to a 1-dimensional
set. We shall do this by constructing a null sequence $U_0,U_1, U_2,"ldots$
of disjoint Peano domains such that, for each $i$,
$p_i "in U_i$, and such that
the union $"cup_iU_i$ is dense in $"St$. Each set 
$"cl(U_i)"setminus "{p_i"}$ can be deformed onto a 1-dimensional
set that contains its boundary by Theorem~2.5.3. 
Since these sets form a null
sequence, the deformations can be combined to find a deformation
that takes $X$ onto the union of $"St "setminus "cup_iU_i$
and the sets $"bd(U_i)$. Each of these sets is a compact 1-dimensional
set. Hence their (countable) union is 1-dimensional.

The domains $U_i$ are created by a long induction. Each step of
the induction constructs a null sequence of Peano domains. 
At step 0 of the induction, an individual domain can have
diameter as large as the diameter of $"St$. Thereafter,
however, we may restrict the maximum diameter of a Peano domain 
at step $i$ to be bounded by $1/i$. Hence the union of this countable
collection of null sequences is also a null sequence.

We consider $"St$
as $"Rt "cup "{"infty"}$. We may assume that $p_0 = "infty "in D(X)$.
By scaling and translating $"Rt$, we find that we may assume that
$[D(X) "setminus "{"infty"}] "cup B(X)$ lies in the interior of
the closed unit square $S = [0,1]"times [0,1]$.

We begin now the construction of our first null sequence of Peano
domains. We outline the strategy. The reader who digests this
strategy will be able to avoid getting lost in the details. We are trying
to fill the open set $U_0 = "St "setminus B(X)$ with 
small Peano domains, more precisely a null sequence of Peano domains,
that are punctured (contain points of $D(X)$). We therefore cover $U_0$
with a fine grid to divide it into small pieces. What happens then
is reminiscent of the children's story, ``Fortunately, Unfortunately.''
Fortunately, some of these small pieces will be punctured.
Unfortunately some will be unpunctured. 
Fortunately, the unpunctured pieces
form a null sequence by hypothesis (2); unfortunately, however, 
they must be attached to adjacent pieces that
are punctured and, unfortunately, the adjacent punctured pieces
need not form a null sequence. Fortunately, we can carve out of 
the adjacent punctured pieces
a null sequence of smaller punctured pieces to which we can
attach the unpunctured pieces. Unfortunately, the process
of carving out small punctured pieces creates new unpunctured
pieces. Fortunately, the new unpunctured pieces form a
null sequence that we can attach to the null sequence
of punctured pieces. Unfortunately, the carving out of small
punctured pieces creates new, as yet unattached, punctured pieces 
that need not form a null sequence. Fortunately, the unattached
punctured pieces are uniformly small and, together, form
a new open set $U_1$ that satisfies hypotheses (1) and (2).
We can then undertake the inductive step with a new open
set whose pieces are smaller than at the previous stage.
Here are the details.

{"bf Step 1. Creating small pieces.} We impose a square grid on $S$
consisting of a large square formed from 
small constituent closed squares.
Since the set $D(X)$ is countable, we lose no generality
in assuming that the edges of the grid miss $D(X)$.
The grid divides the open set $U_0 = "St "setminus B(X)$ into 
many components. We call the collection of such components
$"C_0$. More precisely: {"bf (i)} The set $"St "setminus "nt(S)$
is an element of $"C_0$; {"bf (ii)} If $T$ is any small, closed,
constituent square of the grid, then each component of $T "setminus B(X)$
is also an element of $"C_0$. Note that the elements
of $"C_0$ are not in general disjoint since they can intersect
along the edges of the grid.

{"bf Step 2. Collecting the unpunctured pieces into a
null sequence of small sets.}
 Let $"C_0'$ denote the subcollection of $"C_0$ consisting
of those elements whose interiors are unpunctured. We take the
union $"cup"C_0'$ of the elements of $C_)'$ and claim two things:
{"bf (iii)} The components of $"cup"C_0'$ form a null sequence,
and {"bf (iv)} Each component of $"cup"C_0'$ shares an edge with an
element of $"C_0$ whose interior is punctured.

{"bf Proof of (iii).} We apply here the fundamental
principle of convergence of continua from Section~2.2. The
argument could be repeated almost verbatim perhaps
four more times in the course of Section~2. Often
we will have to consider two cases, depending on whether
the limit continuum contains a point in the interior
of a constituent square of the superimposed grid or does not.
We will not always repeat the details after this first argument.
Here are the details:

Suppose $"e>0$, and suppose that there exist
components $Y_1,Y_2,"ldots$, each of diameter $"ge "e$. We may
assume that $Y_i "to Y$ in the sense of Section~2.2, where $Y$
is a continuum of diameter $"ge "e$. 

Suppose first that $Y$ contains a point in the interior of some
constituent square. Then a small annulus $A$ about that point in
the interior of the constituent square intersects all but
finitely many of the $Y_i$ in a component that crosses $A$
from one boundary component to the other, which easily
gives a contradiction to hypothesis (2).

Suppose next that $Y$ lies in the 1-skeleton of the grid. Then
it contains an interval of an edge of one of the small constituent
squares. In this case, we may take an annulus $A$ that
surrounds an interior point of the interval and intersects
each of the two adjacent squares in a disk (half of an annulus).
Again, all but finitely many of the $Y_i$ will intersect one
of these two disks in a component that crosses the disk
from one side to the opposite, which easily gives a
contradiction to hypothesis (2).

This completes the proof of (iii).

{"bf Proof of (iv).} We may expand the elements of $"C_0$ slightly
without introducing intersections between sets that did not already
intersect; we obtain thus an open covering of $U_0$. Each component
of $U_0$ is punctured, by hypothesis (1). In each component $V$, any
two elements of $"C_0$, as expanded, that lie in $V$ are joined
by a finite chain of such elements by a standard connectedness
argument. A minimal such chain connects each element of $"C_0'$ to
an element of $"C_0$ that is punctured. Property (iv) follows.

{"bf Step 3. Attaching the unpunctured pieces of Step 2 to
a null sequence of punctured pieces.} To each component
$K$ of $"cup "C_0'$ we assign a punctured element $L = L(K) "in "C_0$ 
that intersects $K$ along at least one edge. Such an element
$L(K)$ exists by (iv) of Step~2. The elements $L$ thus chosen
definitely need not form a null sequence, but we shall carve
out from such elements $L$ a new null sequence of punctured
domains to which we may attach the components $K$. Here is
the argument:

For each component $K$, choose an open arc $A(K)$ along which
$K$ is attached to $L(K)$. Choose a point $p(K) "in A(K)$.
Enumerate these points as $q_1,q_2,"ldots$. Each $q_i$ belongs
to a specific $K_i$, and arc $A_i$, and component $L_i = L(K_i)$.

Choose an arc $B_1$ in $L_1$ that joins $q_1$ to $D(X)$ irreducibly.
We may require that $B_1 "cap ("hbox{1-skeleton of grid}) = q_1$ and
that, $"forall$ arcs $B$ having the same properties,
$"diam(B_1) "le 2 "diam(B)$.

Proceed inductively. Choose an arc $B_{k+1}$ in $L_{k+1}$ joining
$q_{k+1}$ to $D(X) "cup B_1 "cup "cdots "cup B_k$ irreducibly.
We may require that $B_{k+1}"cap ("hbox{1-skeleton of grid}) = q_{k+1}$
and that, $"forall$ arcs $B$ having the same properties
$"diam(B_{k+1}) "le 2 "diam(B)$.

We make the following claims about the arcs $B_i$:

{"bf (v)} The arcs $B_1,B_2,"ldots$ form a null sequence.

{"bf (vi)} $"forall "e > 0,","exists k$ such that each component of
$"B(k) = B_{k+1} "cup B_{k+2} "cup "cdots$ has diameter less than
$"e$. 

[Note that (vi) implies (v). Properties (v) and (vi) are
stated separately since (v) is used in the proof of (vi).]

{"bf Proof of (v).} Suppose that (v) is not satisfied. Then there
is a subsequence $B_{i_1},B_{i_2},"ldots$ that converges to
a nondegenerate continuum $B$. [This is our second application
of the fundamental principal of Section~2.2.] We may assume that
the $B_{i_j}$ all lie in the same small constituent square $T$ of
the grid and that their initial endpoints $q_{i_1},q_{i_2},"ldots$
converge to a point $q"in "bd T$. Let $A$ be a small annulus
about $q$ that intersects $T$ in a small disk $A'$. 
All but finitely many of the arcs $B_{i_j}$ cross that disk $A'$
in a large component $B_{i_j}'$. By hypothesis (2), only finitely
many components of $A' "cap U_0$ do not contain a point of $D(X)$.
It follows easily that either some $B_{i_j}'$ is in a component
that contains a point of $D(X)$ or is in a component that contains
another $B_{i_k}'$, with $j > k$. In either case, the diameter
of $B_{i_j}$ can be considerably reduced by shortcutting $B_{i_j}$ to
$D(X)$ or to $B_{i_k}$, a contradiction. This completes the proof of (v).

{"bf Proof of (vi).} We shall make strong use of (v). suppose there
is an $"e>0$ such that each of the sets 
$"B(k) = B_{k+1} "cup B_{k+2} "cup "cdots$ contains a component $Y_k$
of diameter $"ge "e$. We may pick from $Y_k$ a subset $Y_k'$ that
is a finite chain of the arcs $B_1,B_2,"ldots$ and that has diameter
$"ge "e$. Passing to a subsequence if necessary, we may assume that
the sets $Y_k'$ are disjoint and that all lie in the same small
constituent square $T$. If $Y_k' = B_{k_1}"cup "cdots B_{k_"ell}$,
with $k_1 < "cdots < k_"ell$, then we call $q(k) = q_{k_"ell}$ the
initial point of $Y_k'$. We may assume that the initial points
$q(k)$ converge to a point $q"in "bd(T)$. Let $A$ be a small annulus
about $q$ that intersects $T$ in a small disk $A'$. Then each $Y_k'$
is a chain of small arcs crossing $A'$ whose links $B_{k_j}$
all intersect $"bd(T)$. There can be at most two such that are disjoint,
a contradiction. This completes the proof of (vi).

From property (vi) it follows easily that each component $B$ of
$B_1"cup B_2"cup"cdots$ is a tree that lies in a single
small constituent square $T$, contains exactly
one point of $D(X)$, and has, as its leaves, special
attaching points $q_j$ in corresponding attaching arcs
$A_j$ of certain components $K_j$ of $"cup"C_0'$. Furthermore,
these trees $B$ form a null sequence of trees. 
Each component of $"cup"C_0'$ is attached to one of these trees at a leaf. 
We thicken each of these trees slightly and disjointly
so that they still form a null sequence, 
still contain one point of $D(X)$ each,
but now intersect the appropriate attaching arcs $A_j$ in 
neighborhoods $A_j'$ of the attaching points $q_j$. The interiors
of the thickened
trees $B'$ are clearly Peano domains since it is an
easy matter to construct a continuous surjection $f:"Bt "to "cl(B')$
that takes $"nt("Bt)$ homeomorphically onto $"nt(B')$. These
Peano domains will form the cores of the Peano domains that
we are attempting to construct in this stage of the induction.
To them, we must attach the components $K_j$ that we have
described above and also certain sets that we will describe
in the next step.

{"bf Step 4. Attaching the unpunctured components created by removing
the thickened trees of Step~3.} When we remove the thickened
trees $B_i'$ from the components $L = L(K)$, we may create
new components that are unpunctured. We must attach each of
those to an adjacent thickened tree $B_i'$. The following
property establishes the fact that these new unpunctured
domains form a null sequence.

{"bf (vii)} Let $"C_0''$ be the collection of sets defined as follows.
If $K$ is a component of $"cup"C_0'$ and $L = L(K) "in "C_0$,
then the components $M$ of $L "setminus "cup_{i=1}^"infty B_i'$
that contain no points of $"C_0''$ form a null sequence.

{"bf Proof of (vii).} Suppose not. Then there are components 
$M_1,M_2,"ldots$ converging to a nondegenerate continuum $M$.

Suppose first that $M$ has a point $p$ that lies in the
interior of a small constituent square $T$. Since $"cup_iB_i$ is
locally a finite graph away from the edges of the grid, and
a finite graph separates an open set locally into only 
finitely many components, $p "not"in "cup_iB_i$. Hence
there is a small annulus $A$ surrounding $p$ that
contains no point of $"cup_iB_i$. Each $M_i$ crosses $A$ in
a ``large'' set, contained in a component of $A"cap U_0$ 
that contains no points of $D(X)$ and no points of
$"cup_iB_i$. There are only finitely many such, a contradiction.

Suppose finally that $M$ lies in the 1-skeleton
of the grid. Then we may suppose that $M$ contains a nondegenerate
interval $I$ of an edge of a small donstituent square $T$, and
we may assume that each $M_i$ also lies in that square. We may
take a small disk neighborhood $A$ of $I'"subset I$ in $T$ so that
all but finitely many $M_i$ cross $A$ from one side to the other
near $I'$. No point of $"nt(I)$ can lie in $"cup_iB_i$, for
$"cup_iB_i$ separates into
only two components near such a point. Hence, only large $B_i$'s
can be near $I'$. Hence $I'$ has a neighborhood in $A$ missing
$"cup_Bi$. But, by hypothesis (2), all but finitely
many of the components crossing $A$ must contain points of
$D(X)$, a contradiction. 

This completes the proof of (vii).

Each of the components $M$ just discussed share an arc with some
thickened tree $B'$. We attach each component $M$ to such an
adjacent $B'$ along an attaching arc.

{"bf Step 5. Completion of the first null sequence of
Peano domains.} We have at this point created three null
sequences of sets, namely, the components $K$ of $"cup"C_0'$,
the components $B'$ of thickened trees, and the
unpunctured components $M$ that were formed when the
thickened trees were carved out of punctured components of
$"C_0$. Using the attaching arcs described earlier, we can
therefore form a null sequence of domains of the form
$V = "nt(B' "cup K_1"cup K_2 "cup "cdots "cup M_1 "cup M_2 "cup "cdots)$,
where $B'$ is a thickened tree and 
the $K$'s and the $M$'s are attached to $B'$ along attaching arcs.

{"bf (viii)} The sets $V$, which obviously form a null sequence
of sets, are all Peano domains.

{"bf Proof of (viii).} We have already noted that $"nt(B')$
is a Peano domain. Each of the sets $K_i$ is a Peano domain
because, by hypothesis (2) of this theorem, it satisfies
hypothesis (2) of Theorem~2.5.2. We see that the sets $M_j$ are
Peano domains because of the following argument. Suppose
there is a disk $D$ such that the components of $M_j "cap D$
do not form a null sequence. We let $V_1,V_2,"ldots$ denote
a sequence of components converging to a nondegenerate continuum $V$.
We get a contradiction exactly as in the argument for
(vii) above.

We now choose, for the closures of $B'$ and for the closures of
each of the $K_i$'s and each of the $M_j$'s a continuous surjection 
from $"Bt$ as in condition (3) of Theorem~2.5.2. By the proof
of Theorem~2.5.2, as noted in the remark following the statement
of Theorem~2.5.2, we may assume that these maps are 1-1 over
the attaching arcs. It is thus an easy matter to piece these
functions together to get a single continuous surjection from
$"Bt$ onto the closure of 
$V = "nt(B' "cup K_1"cup K_2 "cup "cdots "cup M_1 "cup M_2 "cup "cdots)$
of the kind required by Theorem~2.5.2(3).

This completes the proof of (viii).

{"bf Step 6. Preparing for the next stage of the
induction.} If $K$ is an element of $"C_0$ from
which certain thickened trees $B_i'$ have been
removed, then the remaining punctured components all have
diameter less than or equal to the mesh of the covering
grid. However, they need not form a null sequence. We
simply take the union of the interiors of such elements in $"Rt$
to form a new open set $U_1$. This open set forms the input
to the next stage of the induction. We need to verify
the following fact:

{"bf (ix)} The open set $U_1$ satisfies the two conditions (1) and (2)
with which we began Section~4.

{"bf Proof of (ix).} The remaining components are all subsets of
components of elements of $"C_0$, hence have diameter less
than or equal to the mesh of the covering grid.

Suppose that $D$ is a disk and $D"cap U_1$ has infinitely
many large components $M_i$ that contain no point of $D(X)$.
We may assume $M_i"to M$, $M$ nondegenerate. We argue again
exactly as in the proof of (vii) to obtain a contradiction.

Thus hypothesis (2) is satisfied. Since each component of
$U_1$ is, by hypothesis, punctured, hypothesis (1) is also
satisfied.

This completes the proof of (ix).

{"bf Step 7. The inductive step and the completion of the
proof.} We now recycle the new open set $U_1$ as the set
$U_0$ of the argument just given, but use a grid with much
smaller mesh. We repeat this process inductively, infinitely
often. The completion of the argument is then
clear provided we make the following two remarks:

{"bf (x)} We may require that the point $p_i "in D(X)$ lie in
one of the trees in the $i$th stage of the induction.

Indeed, we may choose the mesh so small that,
if $p_i$ has not been used before stage $i$, then $p_i$ is
the only point of $D(X)$ in a square of the grid and
its neighboring squares, all lying in $U_i$. We can choose
to attach the neighboring squares to the square containing
$p_i$.

{"bf (xi)} Eventually, every point $p$ of 
$"St "setminus (D(X)"cup B(X))$ lies
in the closure of the constructed Peano domains.

Indeed, when squares are sufficiently small, every square
containing $p$ misses $D(X)"cup B(X)$. If $p$ has not already
appeared in the closure of one of our Peano domains, then
$p$ will lie in a component $K$ that contains no point of $D(X)$, hence
will be attached to some thickened tree at that stage.

Thus our proof is complete that we can tile the complement
of $B(X)$ with a null sequence of disjoint Peano domains.
Hence, infinitely many applications of Theorem~2.5.3 show
that $X$ can be deformed by a strong deformation retraction
onto a 1-dimensional set.

"n{\bf 5. Proof of Theorem 1.4.} We are given a codiscrete set $X$.
By Theorem~1.2, $X$ is homotopy equivalent to a planar Peano
continuum $M$. We work with $M$. We must show that the fundamental
group of $M$ embeds in the fundamental group of a 1-dimensional
planar Peano continuum $M'$. 

{"bf The construction of the 1-dimensional planar Peano continuum $M'$.}
We shall associate with $M$ a quotient map $"pi:M "to M'$ onto
a 1-dimensional Peano continuum $M'$ in such a way that each
nondegenerate point preimage $"pi"inv(x)$, for $x "in M'$, is an
arc in $M$ with endpoints in $"bd M$.

{"bf Adjusting $M$.} After slight adjustment, we may assume that $"bd M$
contains no vertical interval.

{"bf Proof.} Suppose $"e > 0$ given. It suffices to
show how to move $M$ a distance $< "e$ such that no homeomorphic
copy of $M$ near the new $M$ can have boundary that
contains a vertical interval as large as $"e$.

Let $[a,b]"times[c,d]$ be a rectangle that contains the $2"e$ 
neighborhood of $M$, with the interval $[a,b]$ horizontal and the
interval $[c,d]$ vertical. Let $c = y_0 < y_1 < "cdots < y_k = d$ be
a partition of $[c,d]$ such that $y_i - y_{i-1}< "e/2$ for each $i$.
Let $[a_j,b_j]"times "{Y_j"}$, $j = 1,"ldots,"ell$, be a collection
of horizontal intervals of length $< "e/2$, all levels
$Y_1,"ldots,Y_"ell,y_0,"ldots,y_k$
being distinct, such that, if $v$ is a vertical interval that joins
adjacent intervals $[a,b]"times "{y_i"}$, then $v$ intersects
some $[a_j,b_j]"times "{Y_j"}$.

Since $"bd M$ is nowhere dense in $"Rt$, there exists a
small open ball $B_j$ in $"Rt"setminus "bd M$ arbitrarily close to
$[a_j,b_j]"times "{Y_j"}$. If such a ball is properly chosen,
it is possible to expand all of the $B_j$'s by an
$"e$-homeomorphism of $"Rt$ so that the image of $B_j$ contains
$[a_j-"e/2,b_j+"e/2]"times"{Y_j"}$. Then every vertical
interval in the image of $M$ that has length $"ge "e$ must contain
a vertical interval $v$ as above, hence must intersect the
image of some $B_j$, hence must intersect $"Rt "setminus "bd M$.

Note that missing the $B_j$'s is an open condition. Hence,
copies of $M$ near this moved $M$ can have no boundary interval
of length $"ge "e$.

{"bf The vertical decomposition of $M$, and the
quotient continuum $M'$.} Let $V$ be a vertical line that
intersects $M$. Let $"G(V)$ denote the set of components of $V "cap M$.
Let $"G = "cup_V"G(V)$. Let $"G_0$ be the trivial
extension of $"G$ to all of $"Rt$. (That is, $"G_0 "setminus "G$
consists of the singleton sets of $"Rt "setminus M$.) Let
$"pi:M "to (M' = M/"G)$ and
$"pi':"Rt "to ("Rt/"G_0)$ be the associated quotient maps.

{"bf Claim 1.} The decomposition $"G_0$ is cellular and upper semicontinuous,
so that $"Rt/"G_0$ is homeomorphic with $"Rt$ by the Moore Decomposition
Theorem 2.3.1. Since each element of $"G$ intersects
$"bd M$, $M' = "pi(M) = "pi'(M)$ is nowhere dense in $"Rt "sim "Rt/"G_0$.
Consequently $M'$ is a 1-dimensional Peano continuum.

{"bf Proof (Claim 1).} Since each element of $"G_0$ is a point or an
arc $"G-0$ is cellular. let $g_1,g_2,"ldots$ be elements of
$"G_0$ containing convergent sequences $x_1,x_2,"ldots "to x$
and $y_1,y_2 "ldots "to y$, with $x_i,y_i "in g_i "in "G_0$.
If $x "ne y$, then we must have $g_i$ a vertical interval
in $M$ for all $i$ sufficiently large. Thus $x$ and $y$ must be
elements of $M$ in the same vertical interval. The vertical intervals 
$g_i$ join $x_i$ to $y_i$. Hence their limits contain a vertical interval
from $x$ to $y$, which must lie in $M$. Thus $x$ and $y$ are
in the same element of $"G_0$, and
$"G_0$ is upper semicontinuous.

The remaining assertions of the claim are easily verified.

{"bf Claim 2.} The projection map $"pi:M "to M'$ induces a map
on fundamental groups that is injective. [This claim is the
central assertion of Theorem~1.4.]

{"bf Proof (Claim 2)}. Let $f:"So "to M$ be a continuous function
such that $f' = ("pi"circ f):"So "to M'$ is nullhomotopic in $M'$
(that is, there is a map $F':"Bt "to M'$ that extends $f'$). 
We must show that $f$ is nullhomotopic in $M$. We may assume 
that $f$ has been standardized in the sense that $f$ restricted
to $f"inv("nt(M))$ is piecewise linear and nowhere vertical. Since
$"bd M$ contains no vertical interval by our previous
adjustment, it follows that $f$ is not vertical on
any subinterval. These adjustments will allow us to give
a rather precise analysis of the maps $f$ and $f'$.

{"bf Analysis of $(f' = "pi "circ f):"So "to M'$.}

{"bf Mapping Analysis Lemma.} Suppose that $f':"So "to M'$
is a nullhomotopic mapping from the circle $"So$ into a 
1-dimensional continuum $M'$. Then there is an upper
semicontinuous decomposition $H$ of $"So$ into compacta
that has the following three properties:

(1) The mapping $f'$ is constant on each element of $H$.

(2) The decomposition $H$ is {"sl noncrossing}. That is,
if $h_1$ and $h_2$ are distinct elements of $H$, then
the convex hulls $"H(h_1)$ and $"H(h_2)$ of $h_1$ and
$h_2$ in the disk $"Bt$ are disjoint. [Equivalently, $h_1$
does not separate $h_2$ on $"So$.]

(3) The decomposition $H$ is {"sl filling}. That is,
the disk $"Bt$ is the union of the convex hulls $"H(h)$
of the elements $h"in H$.

{"bf Proof.} Let $F':"Bt "to M'$ be a map that extends
$f':"So "to M'$. We define 
$$H = "{",h = C "cap "So",|","exists x"in M' "hbox{ such that }C "hbox{ is a component of } F"inv(x)"}.$$

It is obvious that $H$ is an upper semicontinuous decomposition of $"So$
into compacta and that $H$ satisfies property~(1). 

The proof of (2) is easy. If $h_1$ separates $h_2$ on $"So$, and
if $h_1 = C_1 "cap "So$ and $h_2 = C_2 "cap "So$, then
$C_1$ and $C_2$ must intersect, a contradiction.

The proof of (3) will require that we show that $"So/H$ is a
contractible set. Assuming that $"So/H$ is
contractible for the moment, we argue
as follows. Let $H'$ be the collection of sets in $"Rt$ that
are either convex hulls $"H(h)$ of elements of $h"in H$ or are
singleton sets that miss all such convex hulls. Since $H$ is
noncrossing, by (2), it follows easily that $H'$ is a cellular,
upper semicontinuous decomposition of $"Rt$. Let 
$"pi:"Rt "to ("Rt/H'"approx "Rt)$ denote the projection map.
If $H$ were not filling, then the contractible set $"pi("So)
"approx "So/H$ would separate 
the nonempty sets $"pi("Rt "setminus "Bt)$ and $"pi("Bt)"setminus "pi("So)$
in $"Rt/H' "approx "Rt$, a contradiction.

The next paragraphs complete the proof by showing that $"So/H$
is contractible.

Each point $p"in "Bt$ lies on, or in a bounded
complementary domain of, a unique such component $C$ of maximal
diameter. We may redefine $F'$ so that $F'(p) = F'(C)$. This
modification does not alter the decomposition $H$. After
this modification, the nondegenerate components form the
nondegenerate elements of a cellular decomposition $G$ of
$"Rt$; and, by the Moore Decomposition Theorem~2.3.1, the quotient
$"Rt/G$ is homeomorphic with $"Rt$. We denote the
quotient map by $"pi:"Rt "to ("Rt/G \approx "Rt)$.
The modified $F'$ factors through the projection
$"pi|"Bt: "Bt "to "Bt/(G|"Bt)$:
$$F':"Bt "quad"mapright{"pi|"Bt} "quad"Bt/(G|"Bt)"quad "mapright{F''}"quad X.$$

The image $"pi("Bt)$ of the disk $"Bt$ is contractible
because it is a strong deformation retract of the disk
$"pi(2"Bt)"subset "Rt/G.$ [The set $"pi(2"Bt)$ is a disk
since it is a compact set in the plane $"Rt/G$ whose boundary 
is a simple closed curve.]

The image $"pi("Bt)$ of the disk $"Bt$ is 1-dimensional
since (i) it admits the mapping $F'':"pi("Bt)"to M'$ into
a 1-dimensional space $M$ and the
point preimages of $F''$ are totally disconnected, while (ii) a map
that reduces dimension by dimension $k$ must have at least
one point preimage of dimension $k$ [7, Theorem~VI~7].

The images $"pi("Bt)$ and $"pi("So)$ are equal for the following
reasons. Since $"pi("Bt)$ is compact and 1-dimensional, the
open set $"pi("Rt"setminus "Bt)$ is dense in the plane $"Rt/G$.
Hence the image of $"pi("Rt "setminus "nt("Bt))$ is the
entire plane. Consequently, $"pi("So)"supset "pi("Bt)$. The
opposite inclusion is obvious. We conclude that $"pi("So)$
is contractible.

This completes the proof of the Mapping Analysis Lemma.

{"bf Completion of the proof that $f:"So "to M$ is contractible.}

We recall the cellular, upper semicontinuous decomposition 
$"G$ of $"Rt$ that has as its nondegenerate elements the
maximal vertical intervals in $M$ and whose
quotient map $pi:"Rt "to "Rt/"G$ takes $M$ onto $M'$. 
We use the Mapping Analysis Lemma to obtain an
upper semicontinuous decomposition $H$ of $"So$ that
models the shrinking of $f' = ("pi "circ f):"So "to M'$
in the 1-dimensional set $M'$. Since the decomposition
$H$ is noncrossing and filling, we may expand this decomposition
$H$ to a decomposition $G$ of $"Bt$ by taking as elements
the convex hulls in $"Bt$ of the elements of $H$.
The shrinking of $f$ in $M$ will rely on the interplay
between the decompositions $"G$ and $G$.
We shall use the decomposition $G$ of $"Bt$ as a model
on which we shall base the construction of a continuous
function $F:"Bt "to M$ that extends $f:"So "to M$.

If, for each $g"in G$, $f|g"cap "So$ were constant (as is
true for $f'$), we could simply define $F(g) = f(g "cap "So)$.
However, this need not be the case. All that we know is that
$"forall g "in G$, $\exists h(g) "in "G$ such that 
$f(g"cap "So) "subset h(g)$. We need to show how to define
$F|g:g "to h(g)"subset M$ in such a way that the union
$F = "cup"{F|g:g "in G"}$ is a continuous extension of $f$.

If $g$ is a single point, then that point lies in $"So$, and
we may define $F(g) = f(g)$. 

If $g$ is an interval with its
ends in $"So$, then we extend the map $f|"bd g$ linearly
to all of $g$.  

If $g$ is a disk, then we use an {"sl ideal triangulation}
of $g$ in the following way:

The set $g$ is the convex hull $"H(h)$ of
a closed subset $h$ of the unit
circle $"So$. Since $g$ is a disk, $h$ contains at least three points.
An {"sl ideal triangle} is a triangle in $"Bt$ that
has its vertices on $"So$. A collection $"{T_i"}$ of
ideal triangles is said to be an ideal triangulation
of the convex hull $g$ provided that the triangles
have disjoint interiors, have vertices in $h$, and have union
whose intersection with $"nt("Bt)$ is
precisely $g "cap "nt("Bt)$.

{"bf Triangulation Lemma.} If $g = "H(h)$ is a disk, then
$g$ has an ideal triangulation.

{"bf Proof (Triangulation Lemma).} Every point 
$x "in "H(h) "cap "nt("Bt)$ has a neighborhood in $"H(h)$ that
is in the convex hull of a finite collection of points in $h$.
[Hint: every point of a convex hull lies in the hull of a finite
subset; consider separately the case where the point is
in the interior or on the boundary of such a finite polygon.]
Hence, every compact subset of $"H(h)"cap "nt("Bt)$ is in
the convex hull of a finite collection of points in $h$.

Let $C_1 "subset C_2 "subset C_3 "subset "cdots$
be an exhaustion of $"H(h) "cap "nt("Bt)$ by compact sets,
and let $F_1"subset F_2 "subset F_3 "subset "cdots$ be finite subsets
of $h$ such that $C_i "subset "H(F_i)$. It suffices to
show that any ideal triangulation $T_i$ of $"H(F_i)$
extends to an ideal triangulation $T_{i+1}$ of $"H(F_{i+1}$,
for then we may take $T(X) = "cup_{i=1}^"infty T_i$.

To extend $T_i$ to $T_{i+1}$, it suffices to see that we can add
one point $p$ at a time to $F_i$. Since each edge of $T_i$
separates $"Bt$, the domain of $"Bt "setminus |T_i|$ that
contains $p$ is bounded by a single edge $rs$ of $T_i$ followed
by an arc of $"S1$ that contains $p$. We simply add the triangle $prs$
to $T_i$.

This completes the proof of the Triangulation Lemma. 
With the Triangulation Lemma in hand, we are ready to define
$F|g:g "to M$, for the case where $g$ is a disk.

 In this case, we note that $h = g "cap "So$ is a
compact, totally disconnected set having at least three points.
Hence, by the Triangulation Lemma, $g$ has an ideal triangulation
$T(g)$. We define $F$ on $g "cap "So$ to equal $f$. On each
triangle $t_i$ of $T(g)$, we define $F$ to be the linear extension
of $f$ restricted to the three vertices of $t_i$.

{"bf Proof that $F|g$ is continuous for each $g "in G$.} If $F|g$
is not continuous, then $"exists x_1,x_2,"cdots "to x$ in $g$
and $"e > 0$ such that, $"forall i$, $d(F(x_i),F(x))"ge "e$.
Since $F|g"cap"So = f|g"cap "So$ is continuous, we
may assume that each $x_i$ lies in $"nt("Bt)$. Since $F$ is
continuous on any finite union of triangles of $T(g)$ and since
$T(g)$ is locally finite in $"nt("Bt)$, we may assume that
$x"in h = g"cap "So$ and that $x_1,x_2,"ldots$ come from distinct
triangles of $T(g)$. Since these triangles cannot accumulate
at any interior point of $"Bt$, they must, in fact, have diameter
going to $0$ and approach $x$. But then their vertices approach $x$ and,
by linearity, their images approach $x$, a contradiction. Hence
$F|g$ is continuous.

{"bf Proof that $F$ is continuous.} If $F$ is not continuous, then
$"exists",x_1,x_2,"cdots "to x$ in $"Bt$ and $"e > 0$ such that,
$"forall", i$, $d(F(x_i),F(x))"ge "e$. Since, $"forall",g "in G$,
$F|g$ is continuous, we may assume that $x_1,x_2,"ldots,x$ all come from
distinct elements $g_1,g_2,"ldots,g$ of $G$. By continuity of $F|"So = f$, 
we may ignore those $x_i$ in $"So$. Hence, we may assume that
$x_i "in "nt("Bt)$, that $g_i$ is either an arc $t_i$ or a disk,
one of whose triangles $t_i$ contains $x_i$. If the $t_i$
approach $x$, then $x "in "So$, the vertices of the $t_i$ approach $x$, and
the images of the $t_i$ approach $F(x)$ by linearity and the
continuity of $F|"So = f$. Otherwise, we may assume that the
$t_i$ approach an edge $t$ of $g$ that contains $x$. Again,
their vertices approach the vertices of $t$, and the
continuity of $F|"So = f$ and linearity imply that $F(x_i) "to F(x)$,
a contradiction. We conclude that $F$ is continuous.

This completes the proof of Theorem 1.4. We recall the corollary
and question associated with Theorem~1.4:

{"bf Corollary 1.4.1.} If $M$ is a planar Peano continuum, then
the fundamental group of $M$ embeds in an inverse limit of
finitely generated free groups.

{"bf Proof.} This theorem is well-known for 1-dimensional
continua. See, for example, [6] and [2].

{"bf Question 1.4.2.} If $M$ is a planar Peano continuum whose
fundamental group is isomorphic with the fundamental group of
some 1-dimensional planar Peano continuum, is it true that
$M$ is homotopically 1-dimensional?

It is not difficult to see that the projection that we have
given that takes $M$ onto $M'$ does not give a surjection
on fundamental group if $M$ is not homotopically 1-dimensional.
The key issue to resolve here is whether an arbitrary group embedding
into the group of a 1-dimensional continuum 
can always be induced by a continuous map.

We add two final corollaries.

{"bf Corollary 1.4.3.} If $M$ is a planar Peano continuum,
$f:"So "to M$ is a loop in $M$, and $f$ is nullhomotopic
in every neighborhood of $M$ in $"Rt$, then $f$ is
nullhomotopic in $M$. 

{"bf Proof.} It follows easily that $f':"So "to M'$ is
nullhomotopic in each neighborhood of $M'$ in $"Rt$. But
it is well-known [6],[2] that this implies that $f'$ is
nullhomotopic in $M'$. Thus the argument of Theorem~1.4
applies to show that $f$ is nullhomotopic in $M$.

{"bf Corollary 1.4.4.} If $M$ is a planar Peano continuum,
then $"pi_1(M)$ embeds in an inverse limit of free groups.
Those free groups may be taken to be the fundamental
groups of standard neighborhoods (disks with holes) of $M$ in $"Rt$.

{"bf Proof.} Corollary 1.4.4 is a standard corollary to
Corollary 1.4.3.
"v

\centerline{"bf References}


\v
\n\hangindent=1in\hangafter=1 \hbox to 1in{[%
1]\hfil}
J. W. Cannon: 
The recognition problem: what is a topological manifold?. {\sl 
Bull. Amer. Math. Soc.\/} 
84 (
1978), 
832-866.

"v
\n\hangindent=1in\hangafter=1 "hbox to 1in{[%
2]"hfil}
J. W. Cannon, G. R. Conner: 
On the fundamental groups of one-dimensional spaces. {"sl 
"/} 
 (
), 
To appear.

"v
\n\hangindent=1in\hangafter=1 "hbox to 1in{[%
3]"hfil}
J. W. Cannon, G. R. Conner, Andreas Zastrow: 
One-dimensional sets and planar sets are acyclic. In
memory of T. Benny Rushing. {"sl 
Topology Appl."/} 
120 (
2002), 
23-45.

"v
\n\hangindent=1in \hangafter=1 "hbox to 1in{[%
4]"hfil}
John B. Conway: {"sl 
Functions of One Complex Variable II"/}. 
Springer-Verlag, 
New York, Berlin, Heidelberg, 
1995. 

\v
\n\hangindent=1in \hangafter=1 \hbox to 1in{[%
5]\hfil}
Robert J. Daverman: {\sl 
Decompositions of Manifolds\/}. 
Pure and Applied Mathematics Series, V. 124, Academic Press, 
New York, San Francisco, London, 
1986. 

"v
\n\hangindent=1in\hangafter=1 "hbox to 1in{[%
6]"hfil}
M. L. Curtis, M. K. Fort Jr.: 
Homotopy groups of one-dimensional spaces. {"sl 
Proc. Amer. Math. Soc."/} 
8 (
1957), 
577-579.
\

"v
\n\hangindent=1in \hangafter=1 "hbox to 1in{[%
7]"hfil}
Witold Hurewicz and Henry Wallman: {"sl 
Dimension Theory"/}. 
Princeton University Press, 
Princeton, 
1941. 

"v
\n\hangindent=1in\hangafter=1 "hbox to 1in{[%
8]"hfil}
Karimov, Repov"^s, Rosicki, Zastrow: 
On two-dimensional planar compacta which are not homotopically equivalent
to any one-dimensional compactum. {"sl 
preprint."/} 
 (
), 
.

\v
\n\hangindent=1in\hangafter=1 \hbox to 1in{[%
9]\hfil}
R. L. Moore: 
Concerning a set of postulates for plane analysis situs,. {\sl 
Trans. Amer. Math. Soc.\/} 
20 (
1919), 
169-178.

\v
\n\hangindent=1in\hangafter=1 \hbox to 1in{[%
10]\hfil}
R. L. Moore: 
Concerning upper semi-continuous collections of continua. {\sl 
Trans. Amer. Math. Soc.\/} 
27 (
1925), 
416-428.

\v
\n\hangindent=1in \hangafter=1 \hbox to 1in{[%
11]\hfil}
James R. Munkres: {\sl 
Topology, Second Edition\/}. 
Prentice Hall, 
Upper Saddle River, N. J., 
2000. 

\v
\n\hangindent=1in \hangafter=1 \hbox to 1in{[%
12]\hfil}
Edwin H. Spanier: {\sl 
Algebraic Topology\/}. 
Springer-Verlag
New York, Berlin, Heidelberg, 
1989 printing. 
\v
\n\hangindent=1in \hangafter=1 \hbox to 1in{[%
13]\hfil}
Raymond Louis Wilder: {\sl 
Topology of Manifolds\/}. 
American Mathematical Society Colloquium Publications, Volume 32, 
Providence, R.I., 
1963 Edition. 

\par\vfill\bye